\documentclass[centertags,12pt]{amsart}
\usepackage{latexsym}
\usepackage{amssymb}
\usepackage{amscd}
\numberwithin{equation}{section}
\makeatletter
\renewcommand{\subsection}{\@startsection
{subsection}{2}{0mm}{\baselineskip}{-0.25cm}
{\normalfont\normalsize\bf}}
\makeatother
\newtheorem{theorem}{Theorem}[section]
\newtheorem{proposition}[theorem]{Proposition}
\newtheorem{lemma}[theorem]{Lemma}
\newtheorem{corollary}[theorem]{Corollary}
{\theoremstyle{remark}
\newtheorem*{claim*}{Claim}
\newtheorem{remark}[theorem]{Remark}}
\theoremstyle{definition}

\newtheorem*{question*}{Question}
\def\F{\mathbb F}
\def\P{\mathbb P}
\def\N{\mathbb N}

\def\cD{\mathcal D}
\def\cF{\mathcal F}
\def\cH{\mathcal H}
\def\cI{\mathcal I}

\def\cO{\mathcal O}

\def\cX{\mathcal X}

\def\fqss{\mathbb F_{q^6}}
\def\fqs{\mathbb F_{q^2}}
\def\fqc{\mathbb F_{q^3}}
\def\fq{\mathbb F_q}
\def\fqp{\mathbb F_p}

\def\aut{{\rm Aut}}

\def\dim{{\rm dim}}
\def\deg{{\rm deg}}

\begin{document}
\author[Cossidente]{A. Cossidente}
\author[Korchm\'aros]{G. Korchm\'aros}
\author[Torres]{F. Torres}\thanks{1991 Math. Subj. Class.: Primary 11G,
Secondary 14G}\thanks{This research was
carried out within the activity of
GNSAGA of the Italian CNR with the support of the Italian
Ministry for Research and Technology. Part of this paper was
written while Torres was visiting ICTP (Trieste-Italy) supported by
IMPA/Cnpq (Brazil) and ICTP}
\title[Curves covered by the Hermitian curve]{On curves covered by the
Hermitian curve, II}
\address{Dipartimento di Matematica
Universit\`a della Basilicata
via N. Sauro 85, 85100 Potenza, Italy}
\email{cossidente@unibas.it}
\email{korchmaros@unibas.it}
\address{IMECC-UNICAMP, Cx. P. 6065, Campinas-13083-970-SP, Brazil}
\email{ftorres@ime.unicamp.br}
     \begin{abstract}
We classify, up to isomorphism, maximal curves covered by the
Hermitian curve $\cH$ by a prime degree Galois covering. We also compute
the genus of maximal curves obtained by the quotient of $\cH$ by several
automorphisms groups. Finally we discuss the value for the third largest
genus that a maximal curve can have.
     \end{abstract}
\maketitle
\section{Introduction}\label{s1}
In the current study of algebraic geometry in positive characteristic
there is a growing interest to curves which are defined over a finite
field $\F$ and have many $\F$-rational points. Such curves, especially
{\em $\F$-maximal curves}, play indeed a very important role in Coding
theory \cite[Chaper 4, \S7]{go}, \cite[\S8]{lint-geer}, \cite[\S VII.4]{sti2},
\cite[\S3.2.3]{vladut}, and some further motivation for their
investigation also comes from Number Theory \cite{moreno}, \cite{ste} and
Finite Geometry \cite{h}.  Here, maximality of a (projective geometrically
irreducible non-singular algebraic) curve means that the number of its
$\F$-rational points attains the Hasse-Weil upper bound
$$
q^2+1+2qg\, ,
$$
with $g$ being the genus of the curve and $\F$ the field $\fqs$ of order
$q^2$. The majority of work has focused on determining the spectrum for
the genera of maximal curves, and the relevant known results and open
problems concern those maximal curves which have large genus $g$ with
respect to the order $q^2$ of the underlying field; see for example
\cite{ckt}, \cite{fgt}, \cite{g-sti-x}, \cite{geer-vl1}, \cite{geer-vl2},
and \cite{geer-vl3}. The precise upper
bound limit is known to be $q(q-1)/2$ \cite{ihara}. It seems plausible
that only few maximal curves can have genus not too distant from the upper
limit. However, the problem of finding and classifying such maximal curves
is still open and appears to be rather difficult. For an up-dated
discussion of the state of the art, see Section \ref{s5}. The most well
known maximal curve is the so-called Hermitian curve $\cH$ which is
defined by the equation
     \begin{equation}\label{eq1.1}
Y^qZ+YZ^q-X^{q+1}=0\, .
     \end{equation}
In fact, $\cH$ is the only $\fqs$-maximal curve, up to isomorphism, of
genus $q(q-1)/2$ \cite{r-sti} and has a very large automorphism group
(over $\fqs$), namely $\aut(\cH)\cong PGU(3,\fqs)$.
Moreover, all $\fqs$-quotient curves of
$\cH$ are also $\fqs$-maximal, see \cite[Proposition 6]{lachaud}. So, it
would be of interest to solve the above mentioned open problem for
quotient curves arising from automorphism groups of small order. The key
to the solution is to
classify all quotient curves of $\cH$ arising from automorphisms of prime
order; in other words to give a complete classification of Galois
$\fqs$-coverings $\pi:\cH\to \cX$ of prime degree.  This is actually the
main goal in the present paper.
Theorem \ref{thm2.1} states that only five
types of such coverings exist.  It should be noted that almost all
explicit examples of maximal curves in the literature are quotient curves
of $\cH$. This suggests that
Theorem 2.1 might be an essential step toward the complete solution to the
spectrum problem for maximal curves having large genus. For this reason,
Theorem 2.1 also gives some further useful information, such as an
explicit equation for a plane model of the covered curve.  Actually, all
these curves appear in previous work as special cases of wider families, see
Remark \ref{rem2.1}. Some new information will be given in Sections \ref{s3} and
\ref{s4}, such as curves of types (I), (II), (III) and (IV) in Theorem \ref{thm2.1}
can be obtained from subgroups of prime order of $SL(2,\fq)$ while those of class (V)
from the Singer group $S$ of $\cH$. This gives a motivation for the study of
quotient curves arising from subgroups of $SL(2,\fq)$ or from
subgroups of the normaliser in $S$.
In Section \ref{s3}, a variant of the classical Riemann-Hurwitz formula will be
stated which allows a straightforward computation of the genus of all such
but tame curves. However, it remains open the apparently very
involved problem of determining an explicit (possibly singular) plane
model for each of them.

The present research is a continuation of \cite{ckt} in which quotient curves
arising from the subgroups of the Singer group of $\aut(\cH)$ have been investigated.
Computations concerning the genera of certain quotient curves covered by $\cH$ are also
given in \cite{g-sti-x}. However, their methods and results are quite different from
ours, apart from a very few overlapping, see Remarks \ref{rem2.1}, \ref{rem3.1},
\ref{rem4.1}.

Throughout the paper we use the term of a {\em curve} to denote a
projective geometrically irreducible non-singular algebraic curve defined
over the algebraic clousure $\bar\fqs$ of a finite field $\fqs$, of characteristic $p$,
equipped with the Frobenius morphism over $\fqs$.
\section{Classification of curves prime degree Galois\\
covered by the Hermitian curve}\label{s2}
Our purpose is to prove the following theorem.
    \begin{theorem}\label{thm2.1}
Let $\cH$ be the Hermitian curve defined over $\fqs$ and $d$ a prime
number. For a curve $\cX$ of genus $g$ such that $\pi:\cH\to \cX$
is a Galois $\fqs$-covering of degree $d$, we have either
$d=2\neq p$ or $d=p$ or $d\geq 3$ and $(q^2-1)(q^2-q+1)\equiv 0\pmod{d}$.
Moreover, up to $\fqs$-isomorphism one the following cases occurs:
    \begin{enumerate}
\item[(I)] $d=2\neq p$, $\cX$ is the non-singular model of the plane curve
$$
y^q+y=x^{(q+1)/2}\, , \qquad \text{and\quad $g=\frac{(q-1)^2}{4}$}\, ;
$$
\item[(II)] $d=p$ with $q=p^t$, $\cX$ is the non-singular model
of one of the following plane curves
        \subitem(1)
$$
\sum_{i=1}^{t}y^{q/p^i}+\omega x^{q+1}=0\, , \qquad \text{and\quad
$g=\frac{q}{2}(\frac{q}{p}-1)$}\, ,
$$
where $\omega$ is a fixed element of $\fqs$ such that $\omega^{q-1}=-1$;
        \subitem(2)
        $$
y^q+y=(\sum_{i=1}^t x^{q/p^i})^2\, ,\qquad \text{and\quad
$g=\frac{q(q-1)}{2p}$}\, ,
        $$
provided that $p\ge 3$;
\item[(III)] $d\ge 3$ and $q\equiv 1\pmod{d}$, $\cX$ is the non-singular model
of the plane curve
$$
y^q-yx^{2(q-1)/d}+\omega x^{(q-1)/d}=0\, , \qquad \text{and\quad
$g=\frac{q(q-1)}{2d}$}\, ,
$$
where $\omega$ is a fixed element in $\fqs$ such that $\omega^{q+1}=1$;
\item[(IV)] $d\ge 3$ and $q\equiv -1\pmod{d}$, $\cX$ is the non-singular
model of one of the following plane curves
\subitem(1)
$$
y^q+y=x^{(q+1)/d}\, , \qquad\text{and\quad
$g=\frac{(q-1)}{2}(\frac{(q+1)}{d}-1)$}\, ;
$$
\subitem(2)
$$
x^{(q+1)/d}+x^{2(q+1)/d}+y^{q+1}=0\, , \qquad\text{and \quad
$g=(\frac{(q+1)(q-2)}{2d} + 1)$}\, ;
$$
\item[(V)] $d\ge 3$ and $(q^2-q+1)\equiv 0\pmod{d}$, $\cX$ is the
non-singular model of the plane curve
$$
s(x^{q/d}, y^{1/d}, x^{1/d}y^{q/d})=0\, ,
$$
where $s(X,Y,Z):=\prod_{\beta^d=1}(\beta X+\beta^q Y+Z)$, with $\beta$
being a primitive $d$-th root of unity.
In this case $g=\frac{1}{2}(\frac{(q^2-q+1)}{d}-1)$.
    \end{enumerate}
    \end{theorem}
    \begin{remark}\label{rem2.1} $\fqs$-maximal curves with genera as in
the theorem are known to exist, see \cite[Examples D,E,F]{g-sti},
\cite[Thm. 3.1]{geer-vl1}, \cite[Remark 5.2]{geer-vl2}, \cite{ckt},
\cite[Thm. 5.1, Corollary 4.5, Example 5.10]{g-sti-x}. However the interesting
question of determining all such maximal curves is still open, apart from case (I)
for which uniqueness up to $\fqs$-isomorphism has been already proved \cite[Thm. 3.1]{fgt}.
In this context, Theorem 2.1 states the uniqueness for maximal curves prime degree
Galois covered by the Hermitian curve, and it also provides a plane model
for such curves by an explicit equation. For example, for $d=3$
(or, equivalently $q\equiv 2\pmod{3}$), Theorem \ref{thm2.1}(V) states that
   \begin{align*}
s(X,Y,Z) & =  (X+Y+Z)(\beta X+\beta^2 Y+Z)
(\beta^2 X+\beta Y+Z)\\
{} & =  X^3+Y^3+Z^3-3XYZ\, ,\\
\intertext{so that}
s(x^{q/3},y^{1/3},x^{1/3}y^{q/3}) & =
x^q+y+xy^q-3{(xy)}^{(q+1)/3}\, ,
   \end{align*}
defines a plane model of a $\fqs$-maximal curve of genus $(q+1)(q-2)/6$.
     \end{remark}
The proof of Theorem \ref{thm2.1} uses the well-known isomorphism
$\aut(\cH)\cong PGU(3,\fqs)$ (see e.g. \cite{mitchell}, \cite{hartley},
\cite{hoffer}, \cite{sti1})
and it depends on the classification of
subgroups of $PGU(3,\fqs)$ of prime order, see Proposition \ref{prop2.1}. We recall that
$\aut_{\fqs}(\cH)=\aut_{\bar\fqs}(\cH)$, \cite[p. 101]{han-p}.
To have an appropriate description of the actions of such
subgroups on $\cH$, we also need four more plane models of $\cH$
different from (\ref{eq1.1}), namely:
   \begin{enumerate}
\item[(M1)]\quad $X^{q+1}+Y^{q+1}+Z^{q+1}=0$;
\item[(M2)]\quad $Y^qZ-YZ^q+\omega X^{q+1}=0$, where $\omega$ is a fixed
element of $\fqs$ such that $\omega^{q-1}=-1$;
\item[(M3)]\quad $XY^q-X^qY+\omega Z^{q+1}=0$, where $\omega$ is a fixed
element in $\fqs$ such that $\omega^{q+1}=-1$;
\item[(M4)]\quad $XY^q+YZ^q+ZX^q=0$.
    \end{enumerate}
Note that each of the models (M1), (M2) and (M3) are
$\fqs$-isomorphic to (1.1). The model (M4) is $\fqc$-isomorphic to
(M1), cf. \cite[Prop. 4.6]{ckt}.
    \begin{proposition}\label{prop2.1}
Let $C_d=\langle T_d\rangle$ be a subgroup of
$\aut(\cH)\cong PGU(3,\fqs)$ of prime order $d$.
Then $d$ is as in Theorem \ref{thm2.1} and up to conjugacy:
\begin{enumerate}
   \item[(I)] If $d=2\neq p$ and $\cH$ is defined by (\ref{eq1.1}), then
$$
T_d: (X,Y,Z)\mapsto (-X,Y,Z)\, ;
$$
   \item[(II)] Let $d=p$.
       \subitem(1)
If $\cH$ is defined by (M2), then
$$
T_d: (X,Y,Z)\mapsto (X,Y+Z,Z)\, ;
$$
        \subitem(2)
If $p\ge 3$ and $\cH$ is defined by (\ref{eq1.1}), then
$$
T_d: (X,Y,Z)\mapsto (X+Z,X+Y+Z/2,Z)\, ;
$$
   \item[(III)] If $d\ge 3$ and $q\equiv 1\pmod{d}$ and $\cH$ is defined by
(M3), then
$$
T_d: (X,Y,Z)\mapsto (\alpha X,{\alpha}^{-1}Y,Z)\, ;
$$
where $\alpha$ is a primitive $d$-th  root of unity;
   \item[(IV)] If $d\ge 3$ and $q\equiv -1\pmod{d}$ and $\cH$ is defined by (M1),
then we have two
possibilities, either
   \subitem(1)
$$
T_d: (X,Y,Z)\mapsto (\alpha X,Y,Z)\, ,\quad\text{or}
$$
    \subitem(2)
$$
T_d: (X,Y,Z)\mapsto (\alpha X,\alpha^{-1}Y,Z)\, ,
$$
where $\alpha$ is a primitive $d$-th root of unity;

  \item[(V)] If $d\ge 3$ and $(q^2-q+1)\equiv 0\pmod{d}$ and $\cH$ is defined
by (M4), then
$$
T_d: (X,Y,Z)\mapsto (\alpha X, \alpha^q Y, Z)\, ,
$$
where $\alpha$ is a primitive $d$-th root of unity.
   \end{enumerate}
   \end{proposition}
   \begin{proof} An essential tool in the proof is the classification of
all maximal subgroups of $PGU(3,\fqs)$ given by Mitchell \cite{mitchell},
$q$ odd, and by Hartley \cite{hartley}, $q$ even (see also \cite[Ch.\! V]{kleidman},
\cite{hoffer}). This group has order
$q^3(q^3+1)(q^2-1)$. Hence there is a Sylow $d$-subgroup of
$PGU(3,\fqs)$ for each prime divisor $d$ of $q$, $q+1$, $q-1$, and
$(q^2-q+1)$. Since every subgroup of order $d$ is contained in a
Sylow $d$-subgroup and any two Sylow $d$-subgroups are conjugate,
we may choose a Sylow $d$-subgroup $R_d$ for each $d$, and only
consider those subgroups $C_d$ of order $d$ that are contained in
$R_d$.

(I) In $PGU(3,\fqs)$, $q$ odd, elements of order 2 are pairwise
conjugate, and if $\cH$ is given by (\ref{eq1.1}), then $T_d$ in (I) is
an automorphism of order $2$ in $\aut(\cH)$.

(II) We first show that $G:=\aut(\cH)\equiv PGU(3,\fqs$) has either
one or two conjugacy classes of subgroups of order $p$, according
as $p=2$ or $p\ge 3$. Let $\cH$ have equation \ref{eq1.1}. Then a
Sylow $p$-subgroup $R_p$ of $G$ fixes the point
$Q:=(0:1:0)$ and consists of all automorphisms
$$
(X,Y,Z)\to (X+aZ,Y+a^q X+bZ,Z)\, ,
$$
where $a, b\in \fqs$, and $b^q+b=0$. Since no non-trivial element
in $R_p$ fixes a further point of $\cH$, two elements of $R_p$ are
conjugate in $G$ iff they are in the stabilizer $G_Q$ of $Q$.
By \cite[\S10.12]{huppert}, $G_Q$ is the semidirect product
of $R_p$ with a group $H$ of order $(q-1)$ comprising all automorphisms
$$
(X,Y,Z)\to (uX, u^{q+1}Y, Z)\, ,\qquad u\in \fqs^*\, .
$$
Note that the center $Z(R_p)$ of $R_p$ consists of all automorphisms
$$
(X,Y,Z)\to (X,Y+b,Z)\, ,\qquad b\in \fq\, .
$$
A direct computation shows that $Z(R_p)$ is a full conjugacy class
of elements of order $p$ in $G_Q$. For $p=2$, each element of
order $p$ is in $Z(R_p)$. Hence we may assume that $p\ge 3$. Now
let $T_p$ be as in (II)(2). A straightforward computation shows that
the centralizer of $T_p$ in $G_Q$ has order $q^2$, as it consists
of all automorphisms
\begin{equation*}
(X,Y,Z)\to (X+aZ, Y+a^q X+bZ, Z)\, ,\qquad a,b \in \fq\, .\tag{$*$}
\end{equation*}
Hence , the conjugacy class of $T_p$ comprises $q(q^2-1)$ elements
of $G_Q$. Since $q^3=q+q(q^2-1)$, this proves that each
non-central element of $R_p$ is conjugate to $T_p$ under $G_Q$.
Thus $G$ has exactly two conjugacy classes of elements of order
$p$, provided that $p\ge 3$. Now we are in a position to prove II(1).

For $p=2$, the case $a=0$ and $b=1$ in $(*)$ gives $T_2$ in (II)(1). For $p\ge 3$,
we see that $T_p$ above is isomorphic to the automorphism in (II)(1) as follows. We
change the model (\ref{eq1.1}) into the model (M2) via the automorphism $(X,Y,Z)\to
(\omega^{-1} X, \omega^{-1} Y, Z)$ with a fixed $\omega\in \fqs^2$ such that
$\omega^{q-1}=-1$; then, the automorphism for $a=0$ and $b=\omega^{q-1}$ in
$(*)$ is turned into the automorphism in (II)(1).

(III) Let $d\ge 3$ and $q\equiv 1\pmod{d}$, and
denote by $D$ the largest multiplicative subgroup of $\fq$ of order a
power of $d$. Then the automorphisms $(X,Y,Z)\to (uX, u^{-1}Y, Z)$ with
$u$ ranging in $D$, form a subgroup which turns out to be a Sylow
$d$-subgroup $R_d$ of $PGU(3,\fqs)$. Since $R_d$ is cyclic, it has
only one subgroup of order $d$. On the other hand, $R_d$ contains
the automorphism $T_d$ as given in this case and we are done.

(IV) Let $d\ge 3$ and $q\equiv -1\pmod{d}$. First we consider the case $d=3$.
Define $j$ as the
greatest power of $3$
which divides $q+1$. Then a Sylow $3$-subgroup of $PGU(3,\fqs)$
has order $3^{j+1}$. To determine such a Sylow $3$-subgroup $R_3$
explicitely, we adopt the plane model (M1) for $\cH$. Let us introduce the
following automorphisms of $\cH$:
     \begin{align*}
{\phi}_{u,v} & : (X,Y,Z)\to (uX, vY, Z)\, ,\\
{\psi}_{u,v} & : (X,Y,Z)\to  (Z, uX, vY)\, ,\\
{\tau}_{u,v} & : (X,Y,Z)\to  (vY, uZ, X)\, .
     \end{align*}
where $u,v\in{\fqs}$. If both $u$ and $v$ only range in the
subgroup $M$ of order $3^j$ of $\fqs^*$, then the above automorphisms
form a group of order $3^{j+1}$ which is a Sylow $3$-subgroup
$R_3$ of $PGU(3,\fqs)$. Note that $R_3$ is the semidirect
product of $C_3\times C_{3^j}'$ by $C_3''$, where
$C_{3^j}=\langle\phi_{\epsilon,1}\rangle$, $C_{3^j}'=
\langle\phi_{1,\epsilon}\rangle$, $C_3''=\langle\psi_{1,1}\rangle$, and
$\epsilon$ is a primitive third root of unity in $\fqs$. Moreover, the
elements of order $3$ in $R_3$ are $\phi_{u,v}$ with $u^3=v^3=1$,
$\psi_{u,v}$ and $\tau_{u,v}$, with $u,v\in M$. Now let
$\phi$ be an element of order $3$ in
$C_{3^j}\times C_{3^j}'$. It is straightforward to check that $\phi$
is conjugate either to $\phi_{\epsilon,1}$ or to $\phi_{\epsilon,\epsilon^2}$
under a suitable element $\phi_{u,v}$, $u,v\in \fqs^*$. This shows that each
subgroup of $C_{3^j}\times C_{3^j}'$  of order $3$ is either (IV)(1) or
(IV)(2), up to conjugacy. The next step is to check that
$\langle\psi_{u,v}\rangle$, $u,v\in M$ and $(uv)^{3^{j-1}}=1$, is
conjugate to $\langle\phi_{\epsilon,\epsilon^2}\rangle$. Let $w$ be
an element in $M$ such that $w^^3=(uv)^{-1}$. Then the
points $(w:uw^2:1)$, $(\epsilon w:\epsilon^2uw^2:1)$,
$(\epsilon^2w:\epsilon uw^2:1)$ defined over $\fqs$ are the fixed
points of $\psi_{u,v}$. None of these points lies on $\cH$, and they
are the vertices of a triangle. According to \cite{mitchell} and
\cite{hartley}, $PGU(3,\fqs)$ contains an element
that takes this triangle to the fundamental triangle. Then the
conjugate of $\psi_{u,v}$ under the same element belongs to
$C_{3^j}\times C_{3^j}'$, and thus $\langle\psi_{u,v}\rangle=\langle
\phi_{\epsilon,\epsilon^2}\rangle$ up to conjugacy. For
$(uv)^{3^{j-1}}\neq 1$, it turns out instead that the fixed
points $(w:uw^2:1)$, $(\epsilon w:\epsilon^2 uw^2:1)$, $(\epsilon^2 w:
\epsilon uw^2:1)$ of $\psi_{u,v}$ are not defined over $\fqs$, because
$w^3=(uv)^{-1}$ yields $w$ to be in a cubic extension of $\fqs$.
As a consequence of \cite{mitchell} and \cite{hartley}, we have
then that $\psi_{u,v}$ is conjugate to an element
of order $3$ in a Singer subgroup of order $(q^2-q+1)$ of
$PGU(3,\fqs)$. Hence the subgroups $\langle\psi_{u,v}\rangle$ with
$u,v\in M$ but $(uv)^{3^{j-1}}\neq 1$, are pairwise conjugate
under $PGU(3,\fqs)$, and thus each of them is
conjugate to $\langle T_3\rangle$ as given in (V). Now, by
$\tau_{u^{-1},v^{-1}}=\phi_{u,v}^2$, all the above assertions
hold true when $\psi_{u,v}$ is replaced by $\tau_{u,v}$, and this
completes the proof for $d=3$.  In the case $d>3$,
a Sylow $d$-subgroup $R_d$ of $PGU(3,\fqs)$ is $C_{d}\times C_{d}'$, with
$C_d=\langle\phi_{\alpha,1}\rangle$ and
$C_d'=\langle\phi_{1,\alpha}\rangle$. Thus each subgroup of order
$d$ of $R_d$ is either
$\langle\phi_{\alpha,1}\rangle$ or
$\langle\phi_{\alpha,\alpha^{-1}}\rangle$, up to conjugacy.

(V) Let $d\ge 3$ and $(q^2-q+1)\equiv 0\pmod{d}$. Then a Sylow
$d$-subgroup is a subgroup of a Singer subgroup of $PGU(3,\fqs)$,
and hence it is conjugate to $C_d=\langle T_d\rangle$ as given in
(V), see \cite{ck}.
\end{proof}
Now we are in a position to prove Theorem 2.1.

{\it Proof of Theorem 2.1}. For each of the subgroups $C_d=\langle T_d\rangle$
listed in Proposition \ref{prop2.1}, we will determine a $\fqs$-plane model for the
quotient curve $\cH/C_d$, or equivalently the subfield
$\Sigma':=\fqs(\cH/C_d)=\fqs(x',y')$ of the Hermitian function field over $\fqs$
$\Sigma:=\fqs(\cH)=\fqs(x,y)$. Afterwards we compute the genus of $\Sigma'$.

(I) According to Proposition \ref{prop2.1}(I), we define $\Sigma$ by
$y^q+y-x^{q+1}=0$. By considering $x':=x^2$ and $y':=y$ we have that $\Sigma'$ is the fixed
field of $C_d$ and that $(y')^q+y'=(x')^{(q+1)/2}$. For the value of $g$ see
\cite[Prop. VI.4.1]{sti2}.

(II) (1) By Proposition \ref{prop2.1}(II.1), $\Sigma$ is assume to be
$y^q-y+\omega x^{q+1}=0$, with $\omega^{q-1}=-1$. Setting $x':=x$ and $y':=y^p-y$, we have that
$[\Sigma: \Sigma']=p$ and that $\Sigma'$ is the fixed field
of $C_d$. Moreover,
$$
{\rm Tr}_{\fq /\fqp}(y'):=\sum_{i=1}^{t}(y')^{q/p^i}=y^q-y\, ,\qquad
q=p^t\, .
$$
Hence ${\rm Tr}_{\fq /\fqp}(y')+\omega {(x')}^{q+1}=0$. As the polynomial
$Tr_{\fq /\fqp}(Y)+\omega {(x')}^{q+1}$ is irreducible, we obtain the claimed plane
model for $\Sigma'$. For the value on $g$ one proceeds as in \cite[Prop. VI.4.1]{sti2}.

(2) Here, by Proposition \ref{prop2.1}(II.2), $\Sigma$ is defined as in (I) above. Setting
$x':= x^p-x$ and $y':=y-x^2/2$, then $\Sigma'$ is the fixed field of $C_d$. An easy
computation shows that
$$
(y')^q+y'=-({\rm Tr}_{\fq/\fqp}(x'))^2/2\, ,
$$
and hence we obtain an equation defining $\Sigma'$. For the value of $g$ see
\cite[Prop. VI.4.1]{sti2}.

(III) By Proposition \ref{prop2.1}(III), we define $\Sigma$ by $xy^q-x^qy+\omega=0$ with
$\omega^{q+1}=-1$. Let $x':= x^d$, $y':=xy$ and
$$
f(X,Y):= Y^q-YX^{2(q-1)/d}+\omega X^{(q-1)/d}\, .
$$
Then $f(x',y')=0$, and we claim that $f(X,Y)$ is
irreducible in $\bar\fq[X,Y]$. To prove the claim, assume on the
contrary that $f(X,Y)=\prod_i f_i(X,Y)$, where $f_i(X,Y)$ are
irreducible in $\bar\fq[X,Y]$. From $\prod_i f_i(X^d,XY)
=X^{q+1}(X^{q+1}+Y^{q+1}+1)$ follows that
one of the factors on the left hand, say $f_1(X^d,XY)$ is of type
$X^t$, ($0\leq t\le q+1$). Then
$f_1(X,Y)=X^{t/d}$, and $t\equiv 0\pmod{d}$. But $f(X,Y)$ has no factor
of type $X^k$ ($k>1$), and the claim is proved. Next we show that
$\Sigma'$ is the fixed field of $C_d$. Clearly, $\Sigma'\subseteq
{\rm Fix}_{C_d}(\Sigma)$. Note that ${\rm ord}(C_d)=d$
yields $[\Sigma:\Sigma']\le d$. In fact, each element
of $\Sigma$ turns out to be a linear combination
of $1,x,\dots ,x^{d-1}$ over $\Sigma'$. To prove the latter
claim, choose an element $h(x,y)\in\Sigma$. Then
$h(x,y)=a_0(x)+\dots +a_i(x)y^i+\dots +a_n(x)y^n$,
$a_i(x)\in\fqs (x)$. Clearly,
$h(x,y)=b_0(x)+\dots +b_i(x){y'}^i+\dots +b_n(x){y'}^n$, where
$b_i(x)=a_i(x)x^{-i}$. Thus $h(x,y)=b_0(x)+\dots +b_i(x){y'}^i+
\dots +b_b(x){y'}^n$. It remains to show that $b_i(x)$ is a linear
combination of $1,x,\dots ,x^{d-1}$ over $\Sigma'$.
For $b_i(x)\in\fqs [x]$, then
$b_i(x)=b_0+\dots b_ix^i+\dots +b_sx^s$. Replacing
$x'=x^{dj+k}$, ($0\le k<d$) by $x^k{x'}^j$, we see that
$b_i(x)$ is a linear
combination of $1,x,\dots ,x^{d-1}$ over $\Sigma'$. Finally,
let $a(x)=a_0+\dots +a_{d-1}x^{d-1}$, $b(x)=b_0+\dots b_{d-1}x^{d-1}$,
$a_i,b_i\in\Sigma'$. Then there exists $c(x)=c_0+\dots +c_{d-1}
x^{d-1}$, $c_i\in\Sigma'/\fqs$ such that $a(x)/b(x)=c(x)$, and this
completes the proof. Since $T_d$ has two fixed points on $\cH$,
namely $(1:0:0)$ and $(0:1:0)$, from the Riemann-Hurwitz formula applied
to $\cH\to \cH/C_d$, the genus is equal to $(q^2-q)/2d$.

(IV) Here $\Sigma$ is defined by $x^{q+1}+y^{q+1}+1=0$.

(1) Let $x':=x^d$ and $y'=y$. Then $[\Sigma:\Sigma']=d$ and
$\Sigma'$ is the fixed field of $C_d$. There are exactly $q+1$ totally
ramified points in $\cH\to \cH/C_d$, namely $(0:\eta:1)$, with
$\eta^{q+1}=-1$, and we obtain the claimed value for $g$.

(2) Let $x'=x^d$, $y'=xy$, and
$$
f(X,Y):=X^{(q+1)/d}+X^{2(q+1)/d}+Y^{q+1}\, .
$$
Then $f(x',y')=0$, and we see that $f(X,Y)$ is $\bar\fq$-irreducible
arguing as in the proof of (III). Since $T_d$ has no fixed point on
$\cH$, the Riemann-Hurwitz formula applied to $\cH\to \cH/C_d$ gives
$g=(q^2-q+2d-2)/2d$.

(V) By (M4) $\Sigma$ can be defined (over $\fqc$) by
$xy^q+y+x^q=0$. We claim that
$$
f(X,Y):=s(X^{q/d}, Y^{1/d}, X^{1/d}Y^{q/d})
$$
is indeed a polynomial over $\fqs$. To show this we only
need to show that $d\!\mid \! i$ and $d\! \mid\! j$ for each term
$c_{i,j}X_1^iX_2^{j-i}X_3^{d-j}$ in $s(X_1,X_2,X_3)$. Clearly,
$s(\gamma X_1,\gamma^q X_2,X_3)=s(X_1,X_2,X_3)$ for each
$\gamma^d=1$, and
$s(X_1,X_2,X_3)=s(X_2,X_3,X_1)=s(X_3,X_1,X_2)$. Polynomials
satisfying both of the above properties have been investigated
in \cite{ck}. In our case, $s(X_1,X_2,X_3)$ contains
$X_2^d$. Thus \cite[Lemma 6]{ck} yields
$d-j\equiv (q^2+1)(-i)\pmod{d}$, whence $j\equiv qi\pmod{d}$ follows
by $(q^2-q+1)\equiv 0\pmod{d}$. Equivalently, there are integers
$u_{ij}$ and $v_{ij}$ such that $qi-j=du_{ij}$ and
$(q-1)j+i=dv_{ij}$. Then
$$
s(X^{q/d}, Y^{1/d}, X^{1/d}Y^{q/d})=\sum_{}^{}{c_{i,j}
X^{qi+d-j}Y^{j-i+q(d-j)}}\, ,
$$
%
%
and the claim follows. Now, let $x':= x^d$, $y':=  y^d$,
and $f(X,Y)$ as defined before. Then
$$
f(x', y')=
\prod_{\beta^d=1}(\beta  x^q+\beta^q  y+ x y^q)=
(xy^q+y+x^q)\prod_{\beta^d =1, \beta\neq 1}
(\beta  x^q+\beta^q y + x y^q)\, ,
$$
so that $f(x', y')=0$. Moreover, since the product on the right side has $d$
irreducible factors,
the irreducibility of $f(X,Y)$ can be proved by arguing as in case (III).
As $T_d$ has exactly three fixed points, and they are the only
(totally) ramified points, from the Riemann-Hurwitz formula we obtain
$g=(q^2-q-d+1)/2d$.
\section{The genus of maximal curves arising from\\
tame subgroups of $SL(2,\fq)$}\label{s3}
We have already noticed that the group of automorphism $\aut(\cH)$
of the Hermitian curve
$\cH$ is isomorphic to $PGU(3,\fqs)$. From the classification of
subgroups of $PSU(3,\fqs)$ given in \cite{mitchell}, \cite{hartley} and \cite{hoffer}
it follows that $\aut(\cH)$ contains a subgroup
$\Gamma$ isomorphic to $SL(2,\fq)$; moreover, any two such subgroups are
conjugate in $\aut(\cH)$.
Geometrically, $\Gamma$ is contained in the subgroup of $\aut(\cH)$
that preserves a non-incident point-line pair $(P_0,\ell)$, where $P_0\in
\P^2(\fqs)\setminus \cH$ and $\ell$ is its polar line
with respect to the unitary polarity associated with $\cH$.
In particular, $\ell$ is a $\fqs$-rational line meeting $\cH$
in $(q+1)$ pairwise distinct $\fqs$-rational points.

In this section our aim is to compute the genus of the quotient curve of $\cH$ arising
from each tame subgroup of $\Gamma$, see Proposition \ref{prop3.1} (recall that
an automorphism group is called tame if its order is prime to the
characteristic of the base field).
For this purpose, we need at first to give a suitable description
of the action of subgroups of $\Gamma$ on $\cH$. We will use the plane
model (M3) in \S2.

We define the above point-line pair $(P_0,\ell)$ by choosing
$P_0=(0:0:1)$ and $\ell$ as the
line at infinity: $Z=0$. Then the subgroup of automorphisms of
$\P^2(\bar\fqs)$
preserving both $(P_0,\ell)$ and $\cH$, consists of maps of type
        \begin{equation}\label{eq2.1}
(X,Y,Z)\to (aX+bY, cX+dY,Z)\, ,
        \end{equation}
where
$$
ac^q-a^qc=0,\quad bd^q-b^qd=0,\quad bc^q-a^qd=-1,\quad \text{and}\quad
ad^q-b^qc=1\, .
$$
Those maps with $a,b,c,d\in \fq$ and $ad-bc=1$,
form a subgroup isomorphic to $SL(2,\fq)$. We choose
this subgroup to represent $\Gamma$.

Let $G$ be a subgroup of $\Gamma$. The following lemma shows that the
action of $G$ on the affine points of $\cH$ is semi-regular, i.e. each
point-orbit of affine points of $\cH$ under $G$ has length equal to the
order of $G$.
        \begin{lemma}\label{lemma3.1}
Let $\tau\in \Gamma$ and $P\in \cH$ an affine point such that $\tau(P)=P$. Then $\tau$ is
the identity map.
        \end{lemma}
        \begin{proof}
It follows from (\ref{eq2.1}) and the fact that $\alpha^q=\alpha$ for each
$\alpha\in \fq$.
        \end{proof}
From now on we assume that $G$ is tame and investigate the action of $G$ on
the set $\cI:=\ell\cap \cH$,
consisting of all points $(1:m:0), m\in\fq$, together with $(0:1:0)$. Since $\Gamma$
acts on $\cI$ as $PSL(2,\fq)$ in its natural
$2$-transitive permutation representation on the projective line over
$\fq$, we have actually to
consider $\bar G$ instead of $G$, where $\bar G$ is the image
of $G$ under the canonical epimorphism
$$
\phi: \Gamma\cong SL(2,\fq) \to PSL(2,\fq)\, .
$$
Note that the kernel of $\phi$ is trivial for $p=2$, otherwise it is
the subgroup of order $2$ generated by the automorphism
$$
(X,Y,Z)\mapsto (-X, -Y, Z)\, .
$$
Hence either ${\rm ord}(G)=2{\rm ord}(\bar G)$ or ${\rm ord}(G)={\rm ord}(\bar G)$, and
in the later case ${\rm ord}(\bar G)$ must be odd.

According to the classification of subgroups of $PSL(2,\fq)$
\cite[Haupsatz 8.27]{huppert}, the tame subgroup $\bar G$ is one of the
following groups:
      \begin{enumerate}
\item[(3.1)] Cyclic of order $d$, where $d\mid (q+1)$ for $p=2$, or $d\mid (q+1)/2$
for $p\ge 3$\, ;
\item[(3.2)] Dihedral of order $2d$, where $d\mid (q+1)/2$ for $p\ge 3$\, ;
\item[(3.3)] Cyclic of order $d$ where $d\mid (q-1)$ for $p=2$, or
where $d\mid (q-1)/2$ id $p\ge 3$\, ;
\item[(3.4)] Dihedral of order $2d$, where $d\mid (q-1)/2$ for $p\ge 3$\,  ;
\item[(3.5)] The group $Sym_4$ for $q^2\equiv 1\pmod{16}\, , \quad p\ge 5$ ;
\item[(3.6)] The group $Alt_4\quad$ for  $p\ge 5$ ;
\item[(3.7)] The group $Alt_5$ for $q^2\equiv 1\pmod{5}\, , \quad p\ge 7$ .
      \end{enumerate}
We will use the symbols $C_\ell$, $D_\ell$ to denote the cyclic group of
order $\ell$ and the dihedral group of order $2\ell$, respectively. The possibilities
for the action of $\bar G$ on $\cI$ are listed in cases (3.1)-(3.7) below.

{\bf Case 3.1.} Here $\bar G$ has $(q+1)/d$ orbits each of them having
length $d$.

{\bf Case 3.2.} If $q\equiv 3\pmod{4}$, then no involution fixes a point,
and each orbit has length $2d$. If $q\equiv 1\pmod{4}$, then every involution
has two fixed points. Hence just two orbits have length
$d$ and the remaining $(q+1)/2d -1$ orbits have length $2d$.

{\bf Case 3.3.} Here $\bar G$ has two fixed points and the remaining
$(q-1)/d$ orbits have length $d$.

{\bf Case 3.4.} Here $\bar G$ has an orbit of lenght $2$. If $q\equiv
3\pmod{4}$, then the remaining $(q-1)/2d$ orbits have length $2d$.
If $q\equiv 1\pmod{4}$, then just $2$ orbits have length $d$ and the remaining
$(q-1)/2d-1$ orbits have length $2d$.

{\bf Case 3.5.} For a point $P$ on the projective line over $\fq$, let $S$
be the stabilizer of $P$ under $Sym_4$. We show first that $S$ is
either trivial, or isomorphic to any of the following groups: $
C_2, C_3$, or $C_4$. If $S\not\cong C_3$,
then $S$ contains an involution $\tau$ that fixes a  point $Q\neq P$.
Since $PSL(2,\fq)$ is
$2$-transitive on the projective line over $\fq$, we may
assume that $P$ is the infinite point and $Q=P_0$ is the origin.
Then $\tau$ is given by the permutation $X'=-X, Y'=y,Z'=Z$, so $\tau$ is
uniquely determined. This yields that $S$ cannot be isomorphic to
$Alt_4, D_4$ or $D_2$. From the classification of subgroups of $Sym_4$,
it remains to show that $S$ is not isomorphic to $Sym_3$.
Let $g\in S$ be an element of order $3$, then $g$ is
given by $X'=cX+d, Y'=Y, Z'=Z$, with $c,d\in \fq$ and $c^3=1$. Then
$cgc$ is the permutation $cX-d$ which is different from
$c^{-1}X-c^{-1}d$.
On the other hand, the latter permutation is $g^{-1}$. Hence
$cgc\neq g^{-1}$, and this shows that $S\not\equiv Sym_3$.

Let $S\cong C_4$. Then $q\equiv 1\pmod{4}$ and $S$ is generated by
$X'=cX, Y'=Y, Z'=Z$ , with $c^4=1$. Since $Sym_4$ contains exactly three
elements of order $4$, $\bar G$ has just one orbit of length $6$.

Let $S\cong C_3$. Then $q\equiv 1\pmod{3}$ and $S$ can
be assumed to be generated
by $X'=cX, Y'=Y, Z'=Z$, with $c^3=1$. Since $Sym_4$ contains exactly
four subgroups of order $3$, it turns out that $\bar G$ has exactly one
orbit of length $8$.

Let $S\cong C_2$. Then $S$ can be assumed to be generated
by the involution $\tau$ given by $X'=-X, Y'=Y, Z'=Z$. Note that $\tau\in
PSL(2,\fq)$ implies $q\equiv 1\pmod{4}$. Clearly, the orbit of
$P$ under $\bar G$ has length $12$. In particular, the conjucacy
class of $\tau$ has size $6$. Hence $\tau$ is a non-central
involution and $\bar G$ has only one orbit of length $12$.

The above discussion proves the following results:

(I) For $q\equiv 1\pmod{4}$ and $q\equiv 1\pmod{3}$, $\bar G$ has
one orbit of length $6$, one orbit of length $8$,
one orbit of length $12$ and each other orbit has length $24$.

(II) For $q\equiv 1\pmod{4}$ and $q\equiv 2\pmod{3}$, $\bar G$ has
one orbit of length $6$, one orbit of length $12$ and each other orbit
has length $24$.

(III) For $q\equiv 3\pmod{4}$ and $q\equiv 1\pmod{3}$, $\bar G$ has
one orbit of length $8$, and each other orbit has length $24$.

(IV) For $q\equiv 3\pmod{4}$ and $q\equiv 2\pmod{3}$, each orbit
under $\bar G$ has length $24$.

{\bf Case 3.6.} A repetition of the arguments used above
shows that the following cases occur:

(I) For $q\equiv 1\pmod{4}$ and $q\equiv 1\pmod{3}$, $\bar G$ has
two orbits of length $4$, one orbit of length $6$ and
each other orbit has length $12$.

(II) For $q\equiv 1\pmod{4}$ and $q\equiv 2\pmod{3}$, $\bar G$ has
one orbit of length $6$ and each other orbit has length $12$.

(III) For $q\equiv 3\pmod{4}$ and $q\equiv 1\pmod{3}$, $\bar G$ has
two orbits of length $4$ and each other orbit has length $12$.

(IV) For $q\equiv 3\pmod{4}$ and $q\equiv 2\pmod{43}$, each orbit under
$\bar G$ has length $12$.

{\bf Case 3.7.} Similar arguments can be used to prove the following.

(I) For $q\equiv 1\pmod{5}, q\equiv 1\pmod{4}$, and $q\equiv 1\pmod{3}$,
$\bar G$ has one orbit of length $12$, one orbit of length $20$,
one orbit of length $30$ and the remaining orbits have length $60$.

(II) For $q\equiv 1\pmod{5}, q\equiv 1\pmod{4}$, and $q\equiv 2\pmod{3}$,
$\bar G$ has one orbit of length $12$, one orbit of length $30$ and the
remaining orbits have length $60$.

(III) For $q\equiv 1\pmod{5}, q\equiv 3\pmod{4}$, and $q\equiv 1\pmod{3}$,
$\bar G$ has one orbit of length $12$, one orbit of length $30$ and the
remaining orbits have length $60$.

(IV) For $q\equiv 1\pmod{5}, q\equiv 3\pmod{4}$, and $q\equiv 2\pmod{3}$,
$\bar G$ has one orbit of length $12$, and the remaining orbits have length $60$.

(V) For $q\equiv 4\pmod{5}, q\equiv 1\pmod{4}$, and $q\equiv 1\pmod{3}$,
$\bar G$ has one orbit of length $20$, one orbit of length $30$ and the remaining orbits
have length $60$.

(VI) For $q\equiv 4\pmod{5}, q\equiv 1\pmod{4}$, and $q\equiv 2\pmod{3}$,
$\bar G$ has one orbit of length $30$ and the remaining orbits have length $60$.

(VII) For $q\equiv 4\pmod{5}, q\equiv 3\pmod{4}$, and $q\equiv 1\pmod{3}$,
$\bar G$ has one orbit of length $20$, and the remaining  orbits have length $60$.

(VIII) For $q\equiv 4\pmod{5}, q\equiv 3\pmod{4}$, and $q\equiv 2\pmod{3}$, each
orbit under $\bar G$ has length $60$.

Now, the previous case by case analysis of the possible actions of $\bar G$
together with the Riemann-Hurwitz  formula (see Lemma \ref{riemann})
allows us to compute the genus of the quotient
curves associate to $G$, provided that $G$ is tame.
We stress that such curves are $\fqs$-maximal.

To state Lemma \ref{riemann} let $\cX$ denote a curve of genus $g$ and $H$ a subgroup
of $\aut(\cX)$. Let $g'$ be the genus of the quotient curve $\cX/H$ and
suppose that the natural morphism $\pi: \cX\to \cX/H$ is separable. Then
the Riemann-Hurwitz formula applied to $\pi$ states
$$
2g-2=n(2g'-2)+\delta\, ,
$$
where $n$ is the order of $H$ and $\delta$ is the degree of the
ramification divisor $D$ associated to $\pi$. For $P\in \cX$ let
$$
n_P:=\#\{\tau\in H: \tau(P)=P\}\, .
$$
Note that $\#\pi^{-1}(\pi(P))=n/n_P$ and that $n_Q=n_P$ for each $Q\in
\pi^{-1}(\pi(P))$. Now assume that $H$ is tame, so that $p$ does not divide
$n_P$ for each $P\in \cX$, and the multiplicity of $D$ at $P$ is
$(n_P-1)$. As a matter of terminology, the orbit of $P$ is said to be {\em small} if it
consists of less than $n$ elements.
    \begin{lemma}\label{riemann}
If $G$ is a tame subgroup of $\aut(\cX)$ and ${\rm ord}(G)= n$ , then
$$
2g-2=n(2g'-2)+\sum_{i=1}^{s}(n-\ell_i)\, ,
$$
where $\ell_1, \ldots,\ell_s$ are the lenghts of the small orbits of $G$ on $\cX$.
    \end{lemma}
We notice that the above computation generalizes Guerrero's approach
\cite[V.2.5]{far-kra} and it can be deduced from the proof of \cite[V.1.3]{far-kra}.
     \begin{proposition}\label{prop3.1}
Let $G$ denote a tame subgroup of $\Gamma\cong SL(2,\fq)$, $g$ the genus of
the quotient curve $\cH/G$. Then we obtain the following values for $g$, where $\bar
G$ denotes the image of $G$ under the canonical epimorphism $SL(2,\fq)\to PSL(2,\fq)$.
   \begin{enumerate}
\item If $\bar G\cong C_d$, then
$$
g=\begin{cases}
\frac{(q+1)(q-2)}{2d}+1 & \text{if $d$ is odd; $d\mid(q+1)$ for $p=2$, or}\\
{}              &         \text{$d\mid(q+1)/2$ for $p\ge 3$}\, ,\\
\frac{(q+1)(q-3)}{4d}+1 & \text{if $d\mid(q+1)/2$ and $p\ge 3$}\, .
\end{cases}
$$
\item If $\bar G\cong D_d$, with $d\mid(q+1)/2$ and $p\ge 3$, then
$$
g=\begin{cases}
\frac{(q+1)(q-3)}{8d}+1   &  \text{for $q\equiv 3\pmod{3}$}\, ,\\
\frac{(q+1)(q-3)+4d}{8d}  &  \text{for $q\equiv 1\pmod{3}$}\, .
\end{cases}
$$
\item If $\bar G\cong C_d$, then
$$
g=\begin{cases}
\frac{q(q-1)}{2d} &   \text{if $d$ is odd; $d\mid(q-1)$ for $p=2$, or}\\
{}          &          \text{$d\mid (q-1)/2$ for $p\ge 3$}\, ,\\
\frac{(q-1)^2}{4d}     &       \text{if $d\mid (q-1)/2$ and $p\ge 3$}\, .
\end{cases}
$$
\item If $\bar G\cong D_d$, with $d\mid (q-1)/2$ and $p\ge 3$, then
$$
g=\begin{cases}
\frac{(q-1)^2+4d}{8d}   &  \text{for $q\equiv 3\pmod{3}$}\, ,\\
\frac{(q-1)^2}{8d}  &      \text{for $q\equiv 1\pmod{3}$}\, .
\end{cases}
$$
\item If $\bar G\cong Sym_4$, $q^2\equiv 1\pmod{16}$, $p\ge 5$, then
$$
g=\begin{cases}
(q-1)^2/96  & \text{for $q\equiv 1\pmod{4}$ and $q\equiv 1\pmod{3}$}\, ,\\
(q^2-2q+33)/96 & \text{for $q\equiv 1\pmod{4}$ and $q\equiv 2\pmod{3}$}\, ,\\
(q^2-2q+61)/96 & \text{for $q\equiv 3\pmod{4}$ and $q\equiv 1\pmod{3}$}\, ,\\
(q^2-2q+93)/96 & \text{for $q\equiv 3\pmod{4}$ and $q\equiv 2\pmod{3}$}\, .
\end{cases}
$$
\item If $\bar G\cong Alt_4$, $p\ge 3$, then
$$
g=\begin{cases}
(q-1)^2/48     & \text{for $q\equiv 1\pmod{4}$ and $q\equiv 1\pmod{3}$}\, ,\\
(q^2-2q+33)/48 & \text{for $q\equiv 1\pmod{4}$ and $q\equiv 2\pmod{3}$}\, ,\\
(q^2-2q+13)/48  & \text{for $q\equiv 3\pmod{4}$ and $q\equiv 1\pmod{3}$}\, ,\\
(q^2-2q+45)/48 & \text{for $q\equiv 3\pmod{4}$ and $q\equiv 2\pmod{3}$}\, .
\end{cases}
$$
\item If $\bar G\cong Alt_5$ and $q^2\equiv 1\pmod{5}$, $p\ge 7$, then
$$
g=\begin{cases}
(q-1)^2/240  & \text{for $q\equiv 1\pmod{5},\ q\equiv 1\pmod{4}$
and $q\equiv 1\pmod{3}$}\, ,\\
(q^2-2q+81)/240 & \text{for $q\equiv 1\pmod{5},\ q\equiv 1\pmod{4}$
and $q\equiv 2\pmod{3}$}\, ,\\
(q^2-2q+61)/240 & \text{for $q\equiv 1\pmod{5},\ q\equiv 3\pmod{4}$
and $q\equiv 1\pmod{3}$}\, ,\\
(q^2-2q+141)/240 & \text{for $q\equiv 1\pmod{5},\ q\equiv 3\pmod{4}$
and $q\equiv 2\pmod{3}$}\, ,\\
(q^2-2q+97)/240 & \text{for $q\equiv 4\pmod{5},\ q\equiv 1\pmod{4}$
and $q\equiv 1\pmod{3}$}\, ,\\
(q^2-2q+177)/240 & \text{for $q\equiv 4\pmod{5},\ q\equiv 1\pmod{4}$
and $q\equiv 2\pmod{3}$}\, ,\\
(q^2-2q+157)/240 & \text{for $q\equiv 4\pmod{5},\ q\equiv 3\pmod{4}$
and $q\equiv 1\pmod{3}$}\, ,\\
(q^2-2q+237)/240 & \text{for $q\equiv 4\pmod{5},\ q\equiv 3\pmod{4}$
and $q\equiv 2\pmod{3}$}\, .
\end{cases}
$$
    \end{enumerate}
    \end{proposition}
\begin{remark}\label{rem3.1} Comparison with results in \cite{g-sti-x} shows that the
only overlapping concerns Proposition \ref{prop3.1}(1)(4). More precisely, case 1 and
\cite[Example 5.10]{g-sti-x} as well as case 4 and \cite[Example 5.6]{g-sti-x}
coincide. Furthermore, note that $SL(2,\fq)$ is a subgroup of the group $\cF$
introduced in \cite[p. 27]{g-sti-x}; actually $\cF$ is the central product of
$SL(2,\fq)$ with a cyclic group of order $(q+1)$. The results in the present section
give an almost complete answer to the question posed in loc. cit.
\end{remark}
\section{The genus of maximal curves arising from
weakly tame subgroups\\
of the normaliser of a Singer subgroup in $PSU(3,\fqs)$}\label{s4}
The automorphism group $\aut(\cH)$ of the Hermitian curve $\cH$  contains
cyclic groups of order $(q^2-q+1)$; any two such groups are conjugate in
$\aut(\cH)$, \cite{mitchell}, \cite{hartley}, \cite{hoffer}.
These groups and their subgroups are the so-called {\em Singer subgroups}
of $\aut(\cH)$.
Moreover, the normaliser $N=N(\Delta )$ of a Singer subgroup $\Delta$ of order
$(q^2-q+1)$ is a group of order $3(q^2-q+1)$ which is actually the semidirect product of
$\Delta$ with a subgroup $C_3$ of order 3.
Let $\cH$ be given by (M4) (cf. \S2). According to \cite[\S3]{ck}, $\Delta$ can be
chosen as the subgroup generated by
$$
h: (X,Y,Z)\to (\alpha X, \alpha^q, Z)
$$
with $\alpha\in \fqss$ a primitive $(q^2-q+1)$-th root of unity,
while $C_3$ is generated by $(X,Y,Z)\to (Y,Z,X)$.
By \cite[Ch. 4]{short}, the subgroups of $N$ up to conjugacy in $N$ are as follows,
where for $i=0,1,2$, we let $h_i$ denote the automorphism
$$
(X,Y,Z)\to (\epsilon^i Y,\epsilon^{2i} Z,X)
$$
of $\cH$, $\epsilon$ being a primitive third root of unity.
\begin{lemma}\label{lemma4.1}
   \begin{enumerate}
\item[(I)] For every divisor $n$ of $(q^2-q+1)$, the cyclic subgroup $C_n$ of order
$n$, with $C_n=\langle h^{(q^2-q+1)/n}\rangle$;
\item[(II)]
\subitem(1) Let $q\equiv 2\pmod{3}$ and $n\equiv 0\pmod{3}$, or $q\equiv 1\pmod{3}$. For every
divisor $n$ of $(q^2-q+1)$, the subgroup of order $3n$ which is the semidirect product
of $C_n=\langle h^{(q^2-q+1)/n}\rangle$ with $\langle h_0\rangle$.
\subitem(2) Let $q\equiv 2\pmod{3}$ and $n\not\equiv 0\pmod{3}$. For every divisor $n$
of $(q^2-q+1)$, the subgroup $G_i$ ($i=0,1,2$) of order $3n$ which is the
semidirect product of $C_n=\langle h^{(q^2-q+1)/n}\rangle$ with $\langle h_i\rangle$.
   \end{enumerate}
\end{lemma}
The genera of the quotient curves arising from the above subgroups of $\aut(\cH)$ are
given in the following
    \begin{proposition}\label{prop4.1} For any integer $n\ge 3$ satisfying
    $(q^2-q+1)\equiv 0 \pmod{n}$, the quotient curves of the Hermitian curve
    over $\fqs$ arising from the tame subgroups
in the normaliser of the Singer subgroup of $\aut(\cH)$ have the following genera
    \begin{enumerate}
\item $g=((q^2-q+1)/n-1)/2$;
\item $g=(q^2-q+1-n)/6n$ for $q\equiv 2\pmod{3}$ and $n\equiv 0\pmod{3}$ or
$q\equiv 1\pmod{3}$;
\item $g=(q^2-q+1-3n)/6n$ for $q\equiv 2\pmod{3}$ and $n\not\equiv 0\pmod{3}$.
     \end{enumerate}

\end{proposition}
    \begin{proof}
In order to apply the Riemann-Hurwitz formula as stated in Lemma \ref{riemann}, we
take a subgroup $G$ from the list in Lemma \ref{lemma4.1}, and determine its
small-orbits on $\cH$. As (I) was investigated in previous work, see remark below,
we limit ourselves to case (II). Then $G$ has a short orbit $\cO$ of length $3$ consisting of the fixed
points of $h$ which are $(1:0:0)$, $(0:1:0)$ and $(0:0:1)$. For $q\equiv 1\pmod{3}$,
$h_0$ has two fixed points $E_1=(\epsilon:\epsilon^2:1)$ and
$E_2=(\epsilon^2:\epsilon:1)$ on $\cH$. If they belonged to the same orbit under
$G$, then $C_n$ would contain an element that sends $E_1$ to $E_2$, and hence
$\epsilon^n=1$ would follow. On the other hand, $q^2-q+1\equiv 0\pmod{3}$ together
with $q\equiv 1\pmod{3}$ implies $n\not\equiv 0\pmod{3}$. This contradiction shows
that $G$ has further two orbits, ${\cO}':=\{(\beta\epsilon : \beta^q\epsilon^2:1) \mid
\beta^n =1\}$, and ${\cO}''=\{(\beta\epsilon^2:\beta^q\epsilon:1) \mid \beta^n=1\}$.
Now, from Lemma \ref{riemann}, $q^2-q-2=3n(2g-2)+3n-3+2(3n-n)$, and thus
$g=(q^2-q-n+1)/6n$. For $q\equiv 2\pmod{3}$, the picture is richer. Let $n\equiv
0\pmod{3}$. Then $G$ contains the linear transformation $(X,Y,Z)\to (\epsilon
X,\epsilon^2 Y,Z)$, and therefore $h_i$, ($0\le i\le 2$), as defined in Lemma
\ref{lemma4.1}, belongs to $G$. A straightforward computation shows that each of these
automorphisms has three fixed points, namely
\begin{align*}
{\rm Fix}(h_0)= &  \{(\epsilon:\epsilon^2:1), (1:\epsilon:\epsilon^2),
(\epsilon^2:1:\epsilon)\}\, ;\\
{\rm Fix}(h_1)= & \{(1:\epsilon^2:1), (1:1:\epsilon^2), (\epsilon^2:1:1)\}\, ;\\
{\rm Fix}(h_2)= & \{(1:\epsilon:1), (1:1:\epsilon), (\epsilon:1:1)\}\, .
\end{align*}
Note that ${\rm Fix}(h_0)$ is disjoint from $\cH$, while both ${\rm Fix}(h_1)$ and
${\rm Fix}(h_2)$ lie
on $\cH$. Also, $h_0$ induces a $3$-cycle on both ${\rm Fix}(h_1)$ and ${\rm Fix}(h_2)$.
It turns
out that $G$ has two more short orbits, both of length $n$. As before, this gives
$g=(q^2-q-n+1)/6n$. Finally, let $n\not\equiv 0\pmod{3}$. As ${\rm Fix}(h_0)$ is disjoint
from $\cH$, $G_0$ has just one short orbit, namely $\cO$. Thus, Lemma \ref{riemann}
gives $g=(q^2-q+3n+1)/6n$. It remains to consider $G_1$ and $G_2$. Note that $h_i$
($0\le i\le 2$) does not belong to $G_j$ for $i\neq j$. This shows that $G_i$,
($i=1,2$) has exactly four short orbits, namely
$\cO$, ${\cO}_i:=\{(\beta\epsilon^{2i}:\beta^q,1): \beta^n=1\}$,
${\cO}_i':=\{(\beta:\beta^q\epsilon^{2i}:1): \beta^n=1\}$,
${\cO}_i''=\{(\beta,\beta^q,\epsilon^{2i}): \beta^n=1\}$.
From Lemma \ref{riemann}, $q^2-q-2=3n(2g-2)+3n-3+3(3n-n)$, and hence
$g=(q^2-q-3n+1)/6n$.
\end{proof}
\begin{remark}\label{rem4.1} Proposition \ref{prop4.1}(1) has been previously stated
in \cite{ckt} and independently in \cite[Thm. 5.1]{g-sti-x}. Instead , Proposition
\ref{prop4.1}(2)(3) provide new genera for $\fqs$-maximal curves.
\end{remark}
\section{On the third largest genus}\label{s5}
The genus $g$ of a $\fqs$-maximal curve $\cX$ satisfies \cite{ihara},
\cite{sti-x}, \cite{ft1}
$$
g\le g_2:=\lfloor\mbox{$\frac{(q-1)^2}{4}$}\rfloor\qquad\text{or}\qquad
g=g_1:=q(q-1)/2\, .
$$
As remarked in \S1, the Hermitian curve $\cH$ is the only $\fqs$-maximal curve (up to
$\fqs$-isomorphism) with genus $g_1$ and hence $\cH$ is the only maximal
curve having genus as large as possible.
The curves defined by the non-singular models of the following plane curves
$$
y^q+y=x^{(q+1)/2}\qquad \text{$q$ odd}\, ,\qquad \text{and}\qquad
\sum_{i=1}^{t}y^{q/2^i}=x^{q+1}\qquad q=2^t\, ,
$$
have genus $(q-1)^2/4$ and $q(q-2)/4$, respectively. This shows that $g_2$ is the
second largest genus for $\fqs$-maximal curves. For $q$ odd, the above curve is the
only $\fqs$-maximal curve (up to $\fqs$-isomorphism) of genus $(q-1)^2/4$. It seems
plausible that uniqueness also holds true for $q$ even but it has been so far proved
under the additional Condition $(*)$ below (see \cite{at}).
Next we look for the third largest genus $g_3$ that $\cX$ can have. Since
the non-singular model of the curve
$$
y^q+y=x^{(q+1)/3}\, ,\qquad q\equiv 2\pmod{3}\, ,
$$
has genus $(q-1)(q-2)/6$, it is reasonable to search $g_3$ in the
interval
    \begin{equation}\label{eq5.1}
]\lceil\mbox{$\frac{(q-1)(q-2)}{6}$}\rceil,
\lfloor\mbox{$\frac{(q-1)^2}{4}$}\rfloor[\, .
    \end{equation}
In fact, according to \cite[Prop. 2.5]{ft2}, for $q$ odd we have
$$
g_3\le (q-1)(q-2)/4\, .
$$
Recall that $\cX$ is equipped with an
$\fqs$-intrinsic linear series $\cD_\cX$ \cite[\S1]{fgt}.
We have $\dim(\cD_\cX)\ge 2$, equality holding iff $\cX$ is
$\fqs$-isomorphic to the Hermitian curve \cite[Thm. 2.4]{ft2}. Now if the
genus belongs to
(\ref{eq5.1}), then $\dim(\cD_\cX)=3$ \cite[Lemma 3.1]{ckt}. So we look
for $g_3$ among $\fqs$-maximal curves $\cX$ such that $\dim(\cD_\cX)=3$.
In this case, the first three positive Weierstrass non-gaps at $P\in \cX$
satisfy \cite[Prop. 1.5(i)]{fgt}
  \begin{equation}\label{eq5.1'}
m_1(P)<m_2(P)\le q<m_3(P)\, .
  \end{equation}
For $P\in \cX(\fqs)$, we have $m_2(P)=q$ and $m_1(P)\ge q/2$ by $2m_1(P)\ge m_2(P)$.
At this point, we invoke Fuhrmann computations \cite[Anhang \S2]{rainer} concerning the
genus of certain
semigroups of type $\langle m,q,q+1\rangle$. Notice that Fuhrmann's results were
summarized in \cite[Lemma 3.4]{ckt}. It follows that $g$ (the genus of
$\cX$) satisfies
    \begin{equation}\label{eq5.2}
g\le \lfloor\mbox{$\frac{q^2-q+4}{6}$}\rfloor\, ,
    \end{equation}
provided that
     \begin{equation}\label{eq5.3}
m_1(P)\not\in \{\lfloor\mbox{$\frac{q+1}{2}$}\rfloor, q-1\}\, ,\qquad P\in
\cX(\fqs)\, .
     \end{equation}
This leads to investigate some consequences of the following implication
    \begin{equation*}
\text{$\forall P\in \cX(\fqs),\quad m_1(P)=q-1\qquad \Rightarrow\qquad
g< \lfloor\mbox{$\frac{q^2-q+4}{6}$}\rfloor$}\, .\tag{$*$}
    \end{equation*}
      \begin{proposition}\label{prop5.1}
If Condition $(*)$ is satisfied, then $g_3=\lfloor\mbox{$\frac{q^2-q+4}{6}$}\rfloor$.
      \end{proposition}
      \begin{remark}\label{rem5.1} If $q\equiv 2\pmod{3}$, the case $d=3$ in Theorem
      \ref{thm2.1}(V) provides a $\fqs$-maximal curve of genus $(q^2-q-2)/6$. For this curve,
$m_1(P)=q-1$ for each $\fqs$-rational point $P$, see \cite[Prop. 6.4]{ckt}.
Therefore Condition $(*)$ above is not trivial.
       \end{remark}
       \begin{corollary}\label{cor5.1} If Condition $(*)$
is satisfied, then $(q^2-q-2)/6$ is the fourth larger genus that a maximal curve can have
for $q\equiv 2\pmod{3}$.
       \end{corollary}
       \begin{proof} It follows from the theorem and the remark.
       \end{proof}
{\em Proof of Proposition \ref{prop5.1}.} We first
notice that $\fqs$-maximal curves having genus
$\lfloor\mbox{$\frac{q^2-q+4}{6}$}\rfloor$ come from the case $d=3$ in
Theorem \ref{thm2.1}.
Now let $\cX$ be a $\fqs$-maximal curve of genus $g$  such that
$\dim(\cD_\cX)=3$. Then (\ref{eq5.2}), (\ref{eq5.3}) together with the
hypothesis allow us to assume
$m_1(P)=\lfloor\mbox{$\frac{q+1}{2}$}\rfloor$. If $q$ is odd, then
$g=(q-1)^2/4$ \cite[Thm. 2.3]{fgt}; otherwise $g=q(q-2)/4$ \cite{at} and
this completes the proof.
     \begin{remark}\label{rem5.2} By Fuhrmann's results (op. cit.),
a  $\fqs$-maximal curve of genus $\lfloor\mbox{$\frac{q^2-q+4}{6}$}\rfloor$ must have
at least a $\fqs$-rational point $P$ such that
$m_1(P)\in \{\lfloor\mbox{$\frac{2q+2}{3}$}\rfloor,
q-2\}$. For $q\equiv 1\pmod{3}$, Proposition \ref{prop5.3}(5) shows that
$(2q+1)/3$ occurs as a non-gap at certain $\fqs$-rational points.
     \end{remark}
Finally, we discuss necessary conditions for the existence of non-trivial
separable $\fqs$-coverings
$$
\pi: \cH\to \cX
$$
from the Hermitian curve $\cH$ to a ($\fqs$-maximal) curve $\cX$.
    \begin{proposition}\label{prop5.2}
Let $g$ denote the genus of $\cX$. If $g>\lfloor\mbox{$\frac{q^2-q+4}{3}$}\rfloor$,
then $\deg(\pi)=2$, $g=\lfloor\mbox{$\frac{(q-1)^4}{4}$}\rfloor$, and one
of the following holds:
     \begin{enumerate}
\item $\cX$ is the non-singular model of
$y^q+y=x^{(q+1)/2}$ provided that $q$ odd
\item $\cX$ is the non-singular model of $\sum_{i=1}^{t}y^{q/2^i}=x^{q+1}$
provided that $q=2^t$.
     \end{enumerate}
     \end{proposition}
     \begin{proof}
     See \cite{at}.
     \end{proof}
     \begin{proposition}\label{prop5.3} Let $g$ denote the genus of $\cX$. If
$g>\lfloor\mbox{$\frac{q^2-q+6}{8}$}\rfloor$ and
$\deg(\pi)>2$ then
     \begin{enumerate}
\item $\deg(\pi)=3$;
\item $\pi$ is unramified iff $q\equiv 2\pmod{3}$ and $g=(q^2-q+4)/6$;
\item If $\pi$ is ramified, then $g\le (q^2-q)/3$.
     \end{enumerate}
Suppose now that $\pi$ is ramified and that
$g>\lfloor\mbox{$\frac{(q-1)(q-2)}{6}$}\rfloor$. Then
     \begin{enumerate}
\item[4.] If $q\equiv 2\pmod{3}$, $q\ge 5$, then $g=(q^2-q-2)/6$ and $\pi$
is (totally)
ramified at 3 points $P_1, P_2, P_3\not\in\fqs(\cH)$. Moreover for each
$i$, $\pi(P_i)\not\in\cX(\fqs)$ and the Weierstrass semigroup at
$\pi(P_i)$ is given by
$$
\{h/3: h\equiv 0\pmod{3}, h\in S\}\, ,
$$
where
$$
S:= \cup_{j=1}^{q-2}[jq-(j-1),jq]\cup\{0,q^2-2q+2, q^2-2q+3,\ldots\}\, .
$$
In particular, $m_1(\pi(P_i))=(2q-1)/3$ and $m_2(\pi(P_i))=q$.
\item[5.] If $q\equiv 1\pmod{3}$, then $g=(q^2-q)/6$ and $\pi$ is
(totally)
ramified at 2 points $P_1, P_2\in \cH(\fqs)$. The Weierstrass semigroup at
$\pi(P_i)\in \cX(\fqs)$ is given by
$$
\langle (2q+1)/3,q,q+1\rangle\, .
$$
In particular, $m_1(\pi(P_i))=(2q+1)/3$.
\item[6.] If $q\equiv 0\pmod{3}$ and $g=(q^2-q)/6$, then $\pi$ is
(totally)
ramified just at 1 point $P_1\in \cH(\fqs)$; moreover $m_1(\pi(P_1))=2q/3$.
\item[7.] If $\pi$ is normal, i.e. if $\pi$ is Galois, then $\cX$ is $\fqs$-isomorphic to
one of the curves of case $d=3$ in Theorem \ref{thm2.1}.
\end{enumerate}
    \end{proposition}
    \begin{proof} $\deg(\pi)=3$ follows from the Riemann-Hurwitz
formula and the hypothesis on $g$ and $\deg(\pi)$. (2) also follows from
Riemann-Hurwitz. To see (3) we can assume that $p\neq 3$ and that $\pi$
has just one (totally) ramified point. Then $(q^2-q-4)\equiv
0\pmod{3}$, a contradiction.

Now let us assume that $\pi$ is ramified at $P\in \cH$ and let $Q:=\pi(P)$.
By the hypothesis on $g$ we have \cite[Lemma 3.1]{ckt}, \cite[Prop.
1.5]{fgt}
      \begin{equation*}
m_1(Q)<m_2(Q)\le q<m_3(Q)\, .\tag{$**$}
      \end{equation*}
On the other hand, the only possibility for the Weierstrass semigroup $H(P)$ at $P$ is
the above semigroup $S$ (whenever $P\not\in\cH(\fqs)$), and
$H(P)=\langle q,q+1\rangle$ (whenever $P\in \cH(\fqs)$ \cite[Thm.
2]{g-vi}.
Notice that if $h\in H(Q)$, then $3h\in H(P)$; the converse
holds for $p\neq 3$ (see e.g. \cite[Proof of Lemma 3.4]{t}).

{\bf Case $q\equiv 2\pmod{3}$.} We claim that $P\not\in\cH(\fqs)$ and
$Q\not\in\cX(\fqs)$. To see this we first suppose that $P\in \cH(\fqs)$.
Then $m_3(Q)=q+1$ and we have 4 elements in $H(P)$ which are congruent to zero
modulo $3$ and are bounded by $3q+3$. This contradicts $(**)$. Now assume that
$Q\in \cX(\fqs)$ so that $m_3(Q)=q+1$. We then have $3q+3\in H(P)$ so that
$3q+3\ge 4q-3$, i.e. $q=5$. In this case $H(Q)=\{0,3,5,6,7,8\ldots\}$ so
that $g=3$. On the other hand, by \cite[Thm. 2.3]{fgt}, $g=4$. This
contradiction completes the proof. Thus by \cite[Lemma 3.4]{t}
$$
g=\#\{\ell\in \N\setminus S: \ell\equiv 0\pmod{3}\}\, ,
$$
so that $g=(q^2-q-2)/6$ as an easy computation shows. Then the
ramification number of $\pi$ is 6 and so it ramifies at three points. The
statement on Weierstrass semigroups follows from \cite[Proof of Lemma
3.4]{t}.

{\bf Case $q\equiv 1\pmod{3}$.} We claim that $P\in\cH (\fqs)$. For $P\in \cH(\fqs)$,
we have indeed just one element in $H(P)$ which is
$\equiv 0\pmod{3}$ and $\le 3q$. This contradicts $(**)$. Now the proof
can be done as in the previous case, except that
$$
\{h/3: h\equiv 0\pmod{3}, h\in\langle q,q+1\rangle\}=\langle
(2q+1)/3,q,q+1\rangle
$$
follows from \cite[\S2]{rainer}.

{\bf Case $q\equiv 0\pmod{3}$.} Due to wild ramifications , the previous argument does
not allow us any more to use the previous argument to compute the
genus as before. For $g=(q^2-q)/3$, the
ramification number is 4. It follows immediately that $\pi$ is ramified
just at one point $P_1\in \cH$. The non-gaps at $P_1$ less than or equal
to $3q$ turn out to be
$$
\text{either}\quad q,2q-1,2q,3q-2,3q-1,3q\quad\text{or}\quad
q,q+1,2q,2q+1,2q+2,3q\, .
$$
Hence $m_1(\pi(P_1))=2q/3$ and $m_2(\pi(P_1))=q$.
\end{proof}
        \begin{remark}\label{rem5.3}
Let $\cX$ be the curve in Theorem \ref{thm2.1}(II)(2) and $P_0$ the unique
($\fqs$) point over $x=\infty$. Note that $\cD_{\cX}=|(q+1)P_0|$. Now,
 it is easy to see that the first $(p+1)$ positive
Weierstrass non-gaps are $2q/p,\ldots ,pq/p, q+1$. This generalizes
Proposition \ref{prop5.3}(6) and shows that $\dim(\cD_{\cX})=p$.
        \end{remark}
        \begin{remark}\label{rem5.4}
Our final remark concerns the open question of determining all $\fqs$-maximal curves
$\cX$ such that $\dim(\cD_{\cX})=3$. To the list of the known examples
given in \cite[\S6]{ckt}, the non-singular model of the curve
$x^{(q+1)/3}+x^{2(q+1)/3}+y^{q+1}=0$, $q\equiv
2\pmod{3}$, (see Theorem \ref{thm2.1}(IV)(2)), has to be added.
        \end{remark}


\begin{thebibliography}{99999999}

\bibitem[AT]{at} M. Abd\'on and F. Torres, {\em On maximal curves in
characteristic two}, in preparation.

\bibitem[CK]{ck} A. Cossidente and G. Korchm\'aros, {\em The algebraic envelope
associated to a complete arc}, Rend. Circ. Mat. Palermo Suppl. 51 Recent Progress in
Geometry, E. Ballico, G. Korchmaros (Eds.), (1998), 9--24.

\bibitem[CKT]{ckt} A. Cossidente, G. Korchm\'aros and F. Torres, {\em On
curves covered by the Hermitian curve}, submitted (AG/9803029).

\bibitem[Far-Kra]{far-kra} H. M. Farkas and I. Kra, {\em Riemann
surfaces}, Grad. Text in Maths. {\bf 71}, second edition, Springer-Verlag,
1992

\bibitem[F]{rainer} R. Fuhrmann, {\em Algebraische Funktionenk\"orper
\"uber endlichen K\"orpern mit maximaler Anzahl rationaler Stellen},
Ph.D. dissertation, Universit\"at GH Essen, Germany, 1995.

\bibitem[FGT]{fgt} R. Fuhrmann, A. Garcia and F. Torres, {\em On maximal
curves}, J. Number Theory {\bf 67}(1) (1997), 29--51.

\bibitem[FT1]{ft1} R. Fuhrmann and F. Torres, {\em The genus of curves
over finite fields with many rational points}, Manuscripta Math.
{\bf 89} (1996), 103--106.

\bibitem[FT2]{ft2} R. Fuhrmann and F. Torres, {\em On Weierstrass points
and optimal curves}, Rend. Circ. Mat. Palermo Suppl. 51 Recent Progress in
Geometry, E. Ballico, G. Korchmaros (Eds.), (1998), 25--46.

\bibitem[G-Sti]{g-sti} A. Garcia, H. Stichtenoth, {\em Algebraic
function fields over finite fields with many rational places},
IEEE Trans. Inf. Theory {\bf 41}(6), (1995), 1548--1563.

\bibitem[G-Sti-X]{g-sti-x} A. Garcia, H. Stichtenoth and C.P. Xing, {\em
On subfields of the Hermitian function field}, preprint (1998).

\bibitem[G-Vi]{g-vi} A. Garcia and P. Viana, {\em Weierstrass points on
certain non-classical curves}, Arch. Math. {\bf 46} (1986), 315--322.

\bibitem[Geer-Vl1]{geer-vl1} G. van der Geer and M. van der Vlugt,
{\em How to construct curves over finite fields with many points},
Arithmetic Geometry, (Cortona 1994), F. catanese Ed., Cambridge University Press,
Cambridge, 169-189, 1997.

\bibitem[Geer-Vl2]{geer-vl2} G. van der Geer and M. van der Vlugt,
{\em Generalized Reed-Muller codes and curves with many points}, Report
W97-22, Mathematical Institute, University of Leiden, The Netherlands, (
alg-geom/9710016).

\bibitem[Geer-Vl3]{geer-vl3} G. van der Geer and M. van der Vlugt,
{\em  Tables of curves with many points}, April 1998, http:\slash\slash www.wins.uva.nl/geer

\bibitem[Go]{go} V.D. Goppa, {\em Geometry and Codes}, Mathematics and its
applications, 24, Kluwer Academic Publishers, Dordrecht, 1988.

\bibitem[Han-P]{han-p} J.P. Hansen and J.P. Pedersen,
{\em Automorphism groups of Ree type, Deligne-Lusztig curves and function fields},
J. reine angew. Math. {\bf 440} (1993), 99--109.

\bibitem[Har]{hartley} R.W. Hartley, {\em Determination of the ternary
collination groups whose coefficients lie in the $GF(2^n)$}, Annals of
Math. {\bf 27} (1926), 140--158.

\bibitem[H]{h} J.W.P. Hirschfeld, {\em Projective Geometries Over Finite
Fields}, second edition, Oxford University Press, Oxford, 1998.

\bibitem[Hoff]{hoffer} A.R. Hoffer, {\em On unitary collineation groups},
J. Algebra {\bf 22} (1972), 211--218.

\bibitem[Hu]{huppert} B. Huppert, {\em Endliche Gruppen I},
Springer-Verlag, Berlin-Heidelberg-New York, 1967.

\bibitem[Ih]{ihara} Y. Ihara, {\em Some remarks on the number of rational
points of algebraic curves over finite fields}, J. Fac. Sci. Tokio {\bf
28} (1981), 721--724.

\bibitem[L-Geer]{lint-geer} J.H. van Lint and G. van der Geer, {\em
Introduction to Coding Theory and Algebraic Geometry},
DMV Seminar Band 12, Birkh\"ause Verlag, Besel, 1988.

\bibitem[Kl]{kleidman} P.B Kleidman, {\em The maximal subgroups of the low-
dimensional classical groups}, PH.D. Thesis, Cambridge 1987.

\bibitem[La]{lachaud} G. Lachaud, {\em Sommes d'Eisenstein et nombre de
points de certaines courbes alg\'ebriques sur les corps finis}, C.R.
Acad. Sci. Paris {\bf 305}, S\'erie I (1987), 729--732.

\bibitem[M]{mitchell} H.H. Mitchell, {\em Determination of the ordinary
and modular ternary linear groups}, Trans. Amer. Math. Soc. {\bf 12}
(1911), 207--242.

\bibitem[Mo]{moreno} C.J. Moreno, {\em Algebraic Curves over Finite
Fields}, Cambridge University Press, Vol. 97, 1991.

\bibitem[R-Sti]{r-sti} H.G. R\"uck and  H. Stichtenoth, {\em A
characterization of Hermitian function fields over finite fields}, J.
reine angew. Math. {\bf 457} (1994), 185--188.

\bibitem[Short]{short}  M.W. Short, {\em The primitive soluble permutation groups of
degree less then $256$}, LNM 1519, Springer-Verlag, 1992.

\bibitem[Ste]{ste} S.A. Stepanov, {\em Arithmetic of Algebraic Curves},
Consultans Bureau, New York and London, 1994.

\bibitem[Sti1]{sti1} H. Stichtenoth, {\em \" Uber die Automorphismengruppe
eines algebraischen Funktionenk\" orpers von Primzahlcharakteristik},
Arch. Math. {\bf 24} (1973), 527--544 and 615--631.

\bibitem[Sti2]{sti2} H. Stichtenoth, {\em Algebraic Function Fields and Codes},
Springer-Verlag Berlin, 1993.

\bibitem[Sti-X]{sti-x} H. Stichtenoth and C.P. Xing, {\em The genus of
maximal function fields}, Manuscripta Math. {\bf 86} (1995), 217--224.

\bibitem[SV]{sv} K.O. St\"ohr and J.F. Voloch, {\em Weierstrass points and
curves over finite fields}, Proc. London Math. Soc. {\bf 52} (1986), 1--19.

\bibitem[T]{t} F. Torres, {\em On certain $N$--sheeted coverings of curves
and numerical semigroups which cannot be realized as Weierstrass
semigroups}, Comm. Algebra {\bf 23}(11) (1995), 4211--4228.

\bibitem[Tsf-Vla]{vladut} M.A. Tsfasman and S.G. Vladut,
{\em Algebraic-Geometric Codes}, Kluwer Academic Publishers,
Dordrechet-Bosotn-London 1991.
%
\end{thebibliography}
\end{document}